\documentclass{article}
\usepackage[%
journal=FoCM,    
lang=american,   
]{ems-journal}

\numberwithin{equation}{section}

\newtheorem{theorem}{Theorem}
\newtheorem{lemma}{Lemma}
\newtheorem{corollary}{Corollary}

\newtheorem{assumption}{Assumption}

\newtheorem{example}{Example}
\newtheorem{remark}{Remark}

\def\ba{\begin{array}}
\def\ea{\end{array}}
\def\argmin{\mathop{\rm argmin}}

\def\Def{\stackrel{\mathrm{def}}{=}}

\def\dom{{\rm dom \,}}
\def\diam{{\rm diam \,}}
\def\proj{{\rm proj }}

\def\beq{\begin{equation}}
\def\eeq{\end{equation}}

\def\R{\mathbb{R}}
\def\E{\mathbb{E}}

\def\la{\langle}
\def\ra{\rangle}

\def\BR{\begin{remark}}
	\def\ER{\end{remark}}
	
\newcommand{\half}{\mbox{${1 \over 2}$}}

\begin{document}

\title{Universal Reduced-Operator Method and High-Order Global Curvature Bounds}
\titlemark{Universal Reduced-Operator Method and High-Order Global Curvature Bounds}

\emsauthor{1}{
	\givenname{Nikita}
	\surname{Doikov}
	\mrid{1407579}
	\zblid{doikov.nikita}
	\orcid{0000-0003-1141-1625}}{N.~Doikov}
\emsauthor{2}{
	\givenname{Yurii}
	\surname{Nesterov}
	\mrid{203226}
	\zblid{nesterov.yurii}
	\orcid{0000-0002-0542-8757}}{Yu.~Nesterov}

\Emsaffil{1}{
	\department{School of Operations Research and Information Engineering (ORIE)}
	\organisation{Cornell University}
	\rorid{05bnh6r87}
	\address{}
	\zip{}
	\city{}
	\country{}
	\affemail{Nikita.Doikov@cornell.edu}}
\Emsaffil{2}{
	\department{1}{}%
	\organisation{1}{Shenzhen Loop Area Institute and School of Data Sciences (Chinese University of Hong Kong, Shenzhen). Professor emeritus at UCLouvain and Corvinus University of Budapest}%
	\rorid{1}{02d5ks197}
	\address{1}{}%
	\zip{1}{}%
	\city{1}{}%
	\country{1}{}%
	\affemail{1}{Yurii.Nesterov@uclouvain.be}
	}

\classification[68Q25, 90C47]{90C25}


\begin{abstract}
In this paper, we develop a new concept of Global Curvature Bound (GCB) for an arbitrary nonlinear operator 
between abstract metric spaces. We use this notion to characterize the global complexity of
high-order algorithms solving composite variational problems, which include
convex minimization and min-max problems.
We develop the new universal Reduced-Operator Method, 
which automatically achieves 
the fastest universal rate within our class,
while our analysis does not need any specific assumptions about smoothness of the target nonlinear operator.
Every step of our universal method of order $p \geq 1$ requires access to the $p$-th order derivative of the operator and the solution of a strictly monotone doubly regularized subproblem.
For $p = 1$, this corresponds to computing the standard Jacobian matrix of the operator
and solving a simple monotone subproblem, which can be handled using different methods of Convex Optimization. All our results are consequences of the new theorem on the quadratic growth of GCB for general nonlinear operators.
\end{abstract}

\maketitle

\section{Introduction}
\label{SectionIntro}

\noindent
{\bf Motivation.}
Numerical optimization methods stem from
the approximation of a target objective by a model that is easy to construct and minimize
at each iteration of an algorithm. 
A prominent example is \textit{gradient descent}, or the \textit{gradient method},
which approximates the target function by its linear model around a certain point and computes the next point as the minimizer of the regularized model by the square of the Euclidean norm.
In order to prove the convergence rate of gradient descent, one usually assumes that 
the gradient of the objective is Lipschitz continuous, which leads to global bounds
on the quality of the linear approximation and provides the method with progress at each iteration, and consequently, a global rate of convergence~\cite{nesterov2018lectures}.

A generalization of this construction leads to higher-order optimization algorithms, which are
based on Taylor approximations of the target objective augmented by appropriate regularization,
and include, in particular cases, Newton's method with cubic regularization and other approaches, assuming that the Hessian is Lipschitz continuous~\cite{nesterov2006cubic,cartis2011adaptive1,grapiglia2017regularized,nesterov2018lectures}. Tensor methods based on the higher derivatives~\cite{birgin2017worst,nesterov2019implementable,cartis2020sharp}
assume that the higher-order derivatives are Lipschitz.
Despite seeming computational burden, it was shown that third-order tensor methods
can be implemented as efficient as the second-order optimization schemes, utilizing the third-order curvature information but avoiding any heavy calculations~\cite{nesterov2019implementable,nesterov2021superfast}.

However, it soon became evident that the assumed one-to-one relationship between
the order of the Taylor approximation and the corresponding Lipschitz continuous derivative is simply incorrect.
For example, gradient methods perform well on all classes with H\"older-continuous gradients, even without knowing the parameters the problem~\cite{nesterov2015universal,nesterov2024primal},
and second-order methods work universally well for problem classes with both second- and third-order 
H\"older derivatives~\cite{grapiglia2017regularized,grapiglia2020tensor,doikov2024super}, as well as for  variations of quasi-self-concordant functions~\cite{doikov2025minimizing,semenov2025gradient}.

Recently, it was shown that universal gradient methods can work far beyond H\"older continuous classes~\cite{nesterov2025universal}.
Moreover, it has become clear that the assumption on H\"older continuity of the gradient is too restrictive: the sum of two simple univariate functions with different H\"oder degrees does not necessarily belong to \textit{any} H\"older class. The successful framework for analyzing truly universal gradient methods, without limiting them to an a priori restrictive parametric assumption, is based on the notion of the \textit{Global Curvature Bound} (GCB) of the objective, developed in~\cite{nesterov2025universal}. 
It appears that the GCB of a function is an intrinsic measure of the quality of its linear approximation, providing  first-order algorithms with parameter-free complexity bounds. At the same time, in each particular case where  a parametric assumption is used, the convergence results obtained from the GCB analysis recover the correct rate of convergence for the algorithm.

In this work, we generalize the notion of GCB beyond first-order algorithms for minimization problems,
introduced in~\cite{nesterov2025universal}. We focus on one of the most general classes of \textit{globally solvable nonlinear problems}, known as Composite Variational Inequalities~\cite{nesterov2023high}. These types of problems include minimization problems as a particular case, and also cover min-max problems and the problems of solving constrained nonlinear equations. To address these general settings, we develop the concept of the global curvature bound in its most abstract form and demonstrate that it preserves several natural and important properties required for the development of efficient numerical algorithms. Further, we use these developments to construct a new Universal Reduced-Operator method of arbitrary order, which extends the classical extra-gradient method~\cite{korpelevich1977extragradient,nemirovski2004prox} by incorporating at least the Jacobian of the target nonlinear vector field,
while attaining universal parameter-free convergence guarantees.

\vspace{1ex}\noindent
{\bf Contents.} 
The rest of the paper is organized as follows.

In section~\ref{SectionCurvature} we develop a concept of Global Curvature Bound (GCB)
for an arbitrary nonlinear operator between abstract metric spaces, that admit a natural notion of geodesics~\cite{bacak2014convex}.
We show that in this general setting, GCB possesses a very important property of a local quadratic growth (Theorem~\ref{TheoremMain}), under a sole condition that the target metric space is of nonpositve curvature. 
Throughout the paper, the main working examples of metric spaces for us will be normed vector spaces,
for which the theory directly applies.

In section~\ref{SectionHigh}, we utilize the notion of GCB for describing the curvature of high-order derivatives (starting from order zero) of a nonlinear operator between Euclidean vector spaces.
We introduce the \textit{integral smoothing} of GCB that allows us to effectively bound the Taylor approximation of the operator. Finally, we prove several basic properties of the new characterisitcs of high-order curvature, such as their monotonicity and convexity, that will be useful for our analysis.

In section~\ref{SectionReduced} we study the Composite Variational Inequality problem. We use the notion of the Reduced Operator to develop a new algorithm~\eqref{AlgorithmMain} for solving such problems. Its computation needs a solution of an auxiliary {\em doubly regularized} subproblem, which seems to be new.
The method is universal, so it does not need any information about the particular target problem class, except the order $p \geq 1$ of the derivatives of the operator to be used in the computations.
Note that the case $p = 0$ corresponds to the existing well-established Extragradient method~\cite{nemirovski2004prox},
which is not covered by our construction.

The most important example of $p = 1$ correspond to using only the Jacobian of the target operator and can be effectively implemented in practice.  We prove global convergence rates of our method in terms of the reduced operator norm, and in terms of the accuracy certificates that bound the corresponding \textit{merit function} for our problem. We describe the complexity of our algorithm in terms of the GCB of high-order derivatives,
which is very general and does not specify any particular bounds on their growth. At the same time,
we show that in all standard cases, such as H\"older continuity of a certain derivative of the operator, our rate matches the optimal one, known for variational problems in the literature~\cite{bullins2022higher,lin2025perseus,jiang2025generalized,adil2022optimal,nesterov2023high}.

In section~\ref{SectionMonotone} we show how to ensure that the subroblem in our algorithm is monotone for \textit{any order} $p$, which is important for possibility of practical implementation of higher-order ($p \geq 2$) methods.

\section{Global Curvature Bound}
\label{SectionCurvature}

\paragraph{Metric Spaces.}
Let $M$ be an abstract metric space with metric $d(\cdot, \cdot)$ that satisfies the standard axioms. We say that a curve $\gamma: [0, 1] \to M$ connecting two points $x = \gamma(0)$ and $y = \gamma(1)$ is \textit{geodesic}, if
\beq \label{Geodesic}
\ba{rcl}
d(\gamma(s), \gamma(t)) & = & |s - t| d( x, y ), \qquad s, t \in [0, 1].
\ea
\eeq
We will always assume that for any $x, y \in M$ there exists at least one geodesic connecting them, i.e., $M$ is a \textit{geodesic metric space}. For convenience, we denote by ${\cal G}(x, y)$ the set of geodesics between given $x, y \in M$, and, by our assumption, ${\cal G}(x, y)$ is always non-empty.

An important class of metric spaces are those that satisfy the following condition.

\begin{assumption} \label{AssumptionGeo}
	For any two geodesics $\gamma$ and $\eta$, the function
	$ t \mapsto  d(\gamma(t), \eta(t))$,  $t \in [0, 1]$ is convex.
\end{assumption}

The metric spaces that satisfy this assumption sometimes referred to as \textit{Busemann spaces} or spaces with \textit{nonpositive curvature} (see, e.g.,~\cite{bacak2014convex} for a detailed discussion and examples).

\begin{example} \label{ExampleVectorSpace}
	Let $M$ be a normed vector space and we set $d(x, y) := \|x - y\|$.
	Then, it is easy to check that  the straight lines are geodesics:
	$$
	\ba{rcl}
	\gamma(t) & = & [x, y]_t \;\; \Def \;\; x + t ( y  - x), \qquad \gamma \;\; \in \;\; {\cal G}(x, y).
	\ea
	$$
	And if $\| \cdot \|^2$ is strictly convex, the geodesics are unique, so ${\cal G}(x, y)$ consists of one element for any $x, y \in M$. Then, if we consider two arbitrary geodesics,
	$$
	\ba{rcl}
	\gamma(t) = [x_1, y_1]_t \quad \text{and} \quad \eta(t) = [x_2, y_2]_t,
	\qquad x_1, y_1, x_2, y_2 \in M,
	\ea
	$$
	the function
	$$
	\ba{rcl}
	t & \mapsto & d(\gamma(t), \eta(t))
	\;\; = \;\; \| x_1 - x_2 + t( y_1 - x_1 + x_2 - y_2 ) \|
	\ea
	$$
	is convex, and hence Assumption~\ref{AssumptionGeo} is satisfied.
\end{example}

The following space is the classical example that appears in the theory of interior-point methods~\cite{nesterov1994interior},
as applied to linear programming problems.

\begin{example} \label{ExampleLP}
	Let $M$ be the positive orthant $\R^n_{++}$ and consider the Riemannian metric induced 
	by the Hessian of the self-concordant barrier $f(x) = -\sum_{i = 1}^n \ln x^{(i)}$.
	It is known (see Section 6.1 in~\cite{nesterov2002riemannian}) that the corresponding geodesic connecting $x, y \in \R^n_{++}$ is
	$$
	\ba{rcl}
	\bigl[ \gamma(t) \bigr]^{(i)} & = & \exp\bigl(  (1 - t) \ln x^{(i)} + t \ln y^{(i)}  \bigr)
	\;\; = \;\; [ x^{(i)} ]^{1 - t} [ y^{(i)} ]^t, \quad 1 \leq i \leq n,	
	\ea
	$$
	and the distance is given by
	$$
	\ba{rcl}
	d(x, y) & = & \Bigl( \, \sum\limits_{i = 1}^n ( \ln y^{(i)} - \ln x^{(i)} )^2 \Bigr)^{1/2}.
	\ea
	$$
	It is easy to check that~\eqref{Geodesic} is satisfied.
	Then, for two geodesics $\gamma \in {\cal G}(x_1, y_1)$ and $\eta \in {\cal G}(x_2, y_2)$, we have
	$$
	\ba{rcl}
	d(\gamma(t), \eta(t)) & = & 
	\biggl[ \,  \sum\limits_{i = 1}^n 
	\Bigl(  (1 - t) \ln \frac{x_2^{(i)}}{x_1^{(i)}} + t \ln \frac{y_2^{(i)}}{y_1^{(i)}}  \Bigr)^2  \biggr]^{1/2},
	\ea
	$$
	which is convex in $t \in [0, 1]$. Hence, Assumption~\ref{AssumptionGeo} is satisfied.
\end{example}

The following example is important for applications in semi-definite programming~\cite{nesterov1994interior}.

\begin{example} \label{ExampleSDP}
	Let $M$ be the open set $S_{++}^n = \{ X \in \R^{n \times n} \; : \; X = X^{\top} \succ 0 \}$ of symmetric positive-definite matrices
	and consider the Riemannian metric induced by the Hessian of the logarithmic barrier $f(X) = - \ln \det(X)$.
	Then, the corresponding geodesic (see Section~6.3 in~\cite{nesterov2002riemannian}) connecting $X, Y \in S^n_{++}$	is
	$$
	\ba{rcl}
	\gamma(t) & = & X^{1/2} ( X^{-1/2} Y X^{-1/2} )^{t} X^{1/2},
	\ea
	$$
	and the distance is given by
	$$
	\ba{rcl}
	d(X, Y) & = & \Bigl[  \, \sum\limits_{i = 1}^n \ln^2 \lambda_i(X^{-1/2}YX^{-1/2})  \, \Bigr]^{1/2},
	\ea
	$$
	where $\lambda_i(\cdot)$ is the $i$th eigenvalue of a symmetric matrix. It is easy to check directly that~\eqref{Geodesic} holds.
	Then, to verify that the distance function $t \mapsto d(\gamma(t), \eta(t)), t \in [0, 1]$ is convex, for any pair of geodesics $\gamma, \eta$, it is enough to check that for any two matrices, $X, Y \in S^n_{++}$, it holds
	\beq \label{SquareRootBound}
	\ba{rcl}
	d(X^{1/2}, Y^{1/2}) & \leq & \frac{1}{2}d(X, Y),
	\ea
	\eeq
	which can be seen as a consequence of Weyl's inequalities~\cite{weyl1949inequalities} between eigenvalues and singular values.
	Indeed, from~\eqref{SquareRootBound}, and noticing that the linear map $X \mapsto A^{-1/2} X A^{-1/2}$ 
	preserves Riemannian distances while mapping $A$ to identity matrix $I$, we ensure that 
	for any geodesic triangle formed by $A, B, C \in S^n_{++}$, we have
	\beq \label{GeoTriangle}
	\ba{rcl}
	d(\gamma_{AB}(1/2), \gamma_{AC}(1/2))
	& = & d(B^{1/2}, C^{1/2})
	\;\; \overset{(\ref{SquareRootBound})}{\leq} \;\;
	\frac{1}{2} d(B, C),
	\ea
	\eeq
	where $\gamma_{AB} \in {\cal G}(A, B)$ and $\gamma_{BC} \in {\cal G}(B, C)$ are 
	the facets of the geodesic triangle.
	Hence, for any four matrices, $A, B, C, D \in S^n_{++}$ and for $\gamma \in {\cal G}(A, B)$, $\eta \in {\cal G}(C, D)$, $\xi \in {\cal G}(B, C)$, we have
	$$
	\ba{rcl}
	d( \gamma(1/2), \eta(1/2) )
	& \!\!\!\!\leq\!\!\!\! & 
	d( \gamma(1/2), \xi(1/2) ) + d(\xi(1/2), \eta(1/2)) \\
	\\
	& \!\!\!\!\overset{(\ref{GeoTriangle})}{\leq}\!\!\!\! &
	\frac{1}{2} \Bigl[ d(A, C) + d(B, D)  \Bigr]
	\;\; = \;\;
	\frac{1}{2} \Bigl[ d(\gamma(0), \eta(0)) + d(\gamma(1), \eta(1))  \Bigr].
	\ea
	$$
	This means that the continuous function $t \mapsto d(\gamma(t), \eta(t))$ is mid-point convex, hence it is convex,
	and Assumption~\ref{AssumptionGeo} is satisfied. See also~\cite{bhatia2009positive} for a detailed treatment of the Riemannian manifold of positive definite matrices.
\end{example}

A general, well-known fact states that any simply connected, complete Riemannian manifold with nonpositive sectional curvature, called a \textit{Hadamard manifold}, satisfies Assumption~\ref{AssumptionGeo}; a classical example is hyperbolic space~\cite{petersen2006riemannian}.

\paragraph{Operator Curvature.}
Now, let $M$ and $N$ be two distinct metric spaces with metrics $d_M$ and $d_N$, respectively. By abuse of notation, we denote both metrics by $d(\cdot, \cdot)$ whenever no ambiguity arises. Consider an arbitrary mapping $A: M \to N$ between these spaces.
We are interested to quantify the curvature of this mapping. 

For any two points $x, y \in M$ we fix a geodesic $\gamma \in {\cal G}(x, y)$ and compare its image $A(\gamma)$ with a geodesic $\xi \in {\cal G}(A(x), A(y))$ in the target space. We define the following \textit{point-wise measure} of the curvature:
\beq \label{DefDelta}
\ba{rcl}
\Delta_A(x, y) & \Def &  \sup\limits_{0 < t < 1} \frac{1}{t (1 - t)} d\bigl( A(\gamma(t)), \, \xi(t) \bigr).
\ea
\eeq
Note that the definition~\eqref{DefDelta} depends on a particular selection of geodesics, and we assume such a choice to be fixed.
The most important case for us is when geodesics are unique (as in all previous Examples~\ref{ExampleVectorSpace}, \ref{ExampleLP}, \ref{ExampleSDP}),
so that $\Delta_A(x, y)$ is uniquely defined.

Using point-wise measure of curvature, we define the \textit{Global Curvature Bound} (GCB) of mapping $A(\cdot)$ as follows:
\beq \label{DefGCB}
\ba{rcl}
\kappa_A(r) & \Def & \sup\limits_{x, y} \Bigl\{   \Delta_A(x, y) \; : \; d(x, y) \leq r  \Bigr\}.
\ea
\eeq
If $D \Def \diam(M) < +\infty$, denote $\Gamma \Def [0, D]$. Otherwise, $\Gamma \Def \R_+ = \{ t \in \R \, : \, t \geq 0\}$.
Note that $\kappa_A(\cdot)$ is a non-decreasing function on $\Gamma$, and $\kappa_A(0) = 0$. Our main assumption is that $\kappa_A(t)$ is well defined for all $t \in \Gamma$.

\begin{example} \label{ExampleNormedCurvature}
	Let $M, N$ be two normed vector spaces. Set $d_M(x, y) = \| x - y\|_M$, $d_N(a, b) = \|a - b\|_N$
	as in Example~\ref{ExampleVectorSpace}. Consider $A: Q \to N$, where $Q \subseteq M$ is a convex domain.
	Note that we can treat $Q$ as its own metric space with induced metric from $M$.
	Then,
	\beq \label{CurvatureVector}
	\ba{rcl}
	\kappa_A(r) & = & \sup\limits_{\substack{x, y \in Q \\ 0 < t < 1}} 
	\Bigl\{ 
	\frac{1}{t (1 - t)}
	\bigl\| 
	\,
	t A(y) + (1 - t)A(x) - A( x + t (y - x) )
	\,
	\bigr\| \; : \; \|y - x\| \leq r
	\Bigr\}.
	\ea
	\eeq
	When $A$ is an affine mapping, the operation $[\cdot, \cdot]_{t}$ of taking the convex combination commutes:
	$[ A(x), A(y)  ]_{t} \equiv A( [ x, y ]_{t} )$, and thus $\kappa_A(r) \equiv 0$.
	Otherwise, the function $\kappa_A(\cdot)$ serves as a certain measure of curvature (how far the mapping is from being affine).
\end{example}

\begin{example} \label{ExampleFunction}
	In the particular case of the previous example, where $N = \R$, we have a function $f: Q \to \R$.
	The corresponding Global Curvature Bound $\kappa_f(\cdot)$
	was introduced in~\cite{nesterov2025universal}
	to justify the universal complexity of gradient methods applied to minimization problems involving $f$.
\end{example}

Despite to its general nature, the exact GCB has one very important property.

\begin{theorem} \label{TheoremMain}
	Let $A: M \to N$ where $M$ and $N$ are two geodesic metric spaces.
	Assume additionally that $N$ satisfy Assumption~\ref{AssumptionGeo}.
	Then, for any $\beta \in [0, 1]$:
	\beq \label{KappaGrowth}
	\ba{rcl}
	\kappa_{A}(\beta r) & \geq & \beta^2 \kappa_{A}(r), \qquad r \in \Gamma.
	\ea
	\eeq
\end{theorem}
\proof
Let $\frac{1}{2} < \beta < 1$. For $r \in \Gamma$ and some small $\epsilon > 0$, let us choose
points $x, y \in M$ and $\alpha \in [0, 1]$ such that $d(x, y) \leq r$ and
\beq \label{MainBoundGeo1}
\ba{rcl}
d( A(\gamma(\alpha)), \xi(\alpha) ) & \geq & 
\alpha(1 - \alpha) \kappa_A(r)  - \alpha(1-\alpha)\epsilon,
\ea
\eeq
where $\gamma \in {\cal G}(x, y)$ is a geodesic connecting $x$ and $y$, and $\xi \in {\cal G}(A(x), A(y))$ is a geodesic in the image space connecting $A(x)$ and $A(y)$. Without loss of generality, we can assume that $0 < \alpha \leq \frac{1}{2}$, so $\frac{\alpha}{\beta} < 1$.

By the definition of GCB, we also have that
\beq \label{MainBoundGeo2}
\ba{rcl}
d(
A(\gamma(\beta)), 
\xi(\beta)
) & \leq & \beta(1 - \beta) \kappa_A(r).
\ea
\eeq

Now, we note that the curve $\eta(t) := \gamma( \beta t )$, $t \in [0, 1]$
is also a geodesic: $\eta \in {\cal G}( x, \gamma(\beta) )$,
and the distance between its end points is bounded as follows:
$$
\ba{rcl}
d( x, \gamma(\beta) ) 
& \overset{(\ref{Geodesic})}{=} &
\beta d(x, y)
\;\; \leq \;\; \beta r.
\ea
$$

Let us denote by $\phi$ the geodesic in the image space connecting the corresponding end points: $\phi \in {\cal G}(A(x), A(\gamma(\beta)))$.
Then, by the definition of GCB and selecting the intermediate point 
$$
\boxed{
	\ba{rcl}
	s & := & \frac{\alpha}{\beta} \;\; \in \;\; (0, 1),
	\ea
}
$$ 
we get that
\beq \label{MainBoundGeo3}
\ba{rcl}
d( A(\gamma(\alpha)), \phi( s ) )
& \!\! = \!\! &
d( A( \eta( s ) ), \phi( s ) )
\leq 
s(1 - s) \kappa_A(\beta r)
= 
\frac{\alpha (\beta - \alpha)}{\beta^2} \kappa_A(\beta r).
\ea
\eeq
Let us observe that the curve $\psi(t) := \xi(\beta t), t \in [0, 1]$ is a geodesic in the image space: $\psi \in {\cal G}( A(x), \xi(\beta) )$. Therefore, we have two geodesics $\psi$ and $\phi$  coming from the same point $\psi(0) = \phi(0) = A(x)$.
We use our Assumption~\ref{AssumptionGeo}. Namely, we know that the function
$$
\ba{rcl}
t \mapsto d( \psi(t), \phi(t) ), \qquad t \in [0, 1],
\ea
$$	
is convex, which gives
\beq \label{ConvGeod}
\ba{rcl}
d( \phi(s), \psi(s))
& \leq &
s d( \phi(1), \psi(1)) 
+ (1 - s) d(\phi(0), \psi(0)) \\
\\
& = &	
s d(\phi(1), \psi(1) ) 
\;\; = \;\;
\frac{\alpha}{\beta}d( A(\gamma(\beta)), \xi(\beta) ) \\
\\
& \overset{(\ref{MainBoundGeo2})}{\leq} & 
\alpha(1 - \beta) \kappa_A(r).
\ea
\eeq
Since, $\psi(s) = \xi(\alpha)$, we finally obtain 
\beq \label{MainBoundGeo4}
\ba{rcl}
d( \phi(s), \xi(\alpha) ) & \leq & 	\alpha(1 - \beta) \kappa_A(r).
\ea
\eeq
Then, using the triangle inequality, we conclude that
$$
\ba{rcl}
\kappa_A(\beta r) & \overset{(\ref{MainBoundGeo3})}{\geq} &
\frac{\beta^2}{\alpha(\beta - \alpha)}
d( A(\gamma(\alpha)), \phi(s) ) \\
\\
& \geq & 
\frac{\beta^2}{\alpha(\beta - \alpha)}
\Bigl[
d( A(\gamma(\alpha)), \xi(\alpha) )
- d( \phi(s), \xi(\alpha) )
\Bigr] \\
\\
& \overset{(\ref{MainBoundGeo1}), (\ref{MainBoundGeo4})}{\geq} &
\frac{\beta^2}{\alpha(\beta - \alpha)}
\Bigl[
\alpha (1 - \alpha) \kappa_A(r) - \alpha(1-\alpha)\epsilon - \alpha(1 - \beta) \kappa_A(r)
\Bigr] \\
\\
& = & 
\beta^2 \kappa_A(r) - \frac{\beta^2 (1 - \alpha) }{\beta - \alpha} \epsilon 
\;\; \geq \;\;
\beta^2 \kappa_A(r)
-
\frac{\beta^2}{(\beta - 1/2)} \epsilon.
\ea
$$
It remains to choose $\epsilon \to 0$.

If $0 < \beta \leq \frac{1}{2}$, then we can represent it as $\beta = \bar{\beta}^k$ with some $\frac{1}{2} < \bar{\beta} < 1$
and $k \geq 1$, and repeat the above reasoning for $\bar{\beta}$ recursively $k$ times.
\qed

We conclude with the following important consequences of our definition that we use in our further developments.

\begin{remark} \label{RemarkLinearBound}
	Assume that $M, N$ are two normed vector spaces as in Example~\ref{ExampleNormedCurvature} and consider $A: \dom A \to N$, with $\dom A \subseteq M$.
	Denote by $DA(x)[h]$ the directional derivative of mapping $A(\cdot)$
	at $x \in \dom A$ along direction $h \in M$. Let us assume that 
	$DA(x)[y - x]$ exists for all $x, y \in \dom A$. Since
	$$
	\ba{rcl}
	\lim\limits_{t \to 0} 
	\frac{1}{t(1 - t)} \bigl[ t A(y) + (1 - t)A(x) - A(x + t(y - x))  \bigr]
	& = & 
	A(y) - A(x) - DA(x)[y - x],
	\ea
	$$
	we have that
	\beq \label{eq-AppUp}
	\ba{rcl}
	\|A(y) - A(x) - DA(x)[y - x] \| & \leq & \kappa_A( \|y - x\| ), \qquad x, y \in \dom A.
	\ea
	\eeq
\end{remark}

We see by~\eqref{eq-AppUp} that the exact GCB ultimately bounds the linear approximation of our operator $A$,
when working in normed vector spaces. Note that we did not assume any bounds on the global behaviour of $A$.
At the same time, when a certain global smoothness assumption holds, we can use this information to specify
the growth rate of GCB, as in the next example.

\begin{example} \label{ExampleHolder}
	As in Remark~\ref{RemarkLinearBound}, let $M$ and $N$ be normed vector spaces,
	and assume that the operator $A: \dom A \to N$, with $\dom A \subseteq M$ is differentiable
	and its derivative is H\"older continuous of some degree $\nu \in [0, 1]$ with constant $H_{\nu} > 0$:
	\beq \label{HolderDA}
	\ba{rcl}
	\max\limits_{h \in M \; : \; \| h \| \leq 1}
	\| DA(x)[h] - DA(y)[h] \|
	& \leq & H_{\nu} \|x - y\|^{\nu}, \qquad x, y \in \dom A.
	\ea
	\eeq
	Then
	\beq \label{HolderGCB}
	\ba{rcl}
	\kappa_A(r) & \leq & \frac{2^{1 - \nu} H_{\nu}}{1 + \nu} r^{1 + \nu}.
	\ea
	\eeq
\end{example}
\proof
By inequality~\eqref{HolderDA}, we have, for any $x, y \in \dom A$:
\beq \label{HolderABound}
\ba{rcl}
\| A(y) - A(x) - DA(x)[y - x] \| & \leq & \frac{H_{\nu}}{1 + \nu} \|y - x\|^{1 + \nu}.
\ea
\eeq
Let us fix two points $x, y \in \dom A$, and apply inequality~\eqref{HolderABound} for each of them and 
the point $[x, y]_t := x + t(y - x)$,
for $t \in (0, 1)$.
We get,
\beq \label{TwoHolderBounds}
\ba{rcl}
\| A(y) - A([x, y]_t) - (1 - t) DA([x, y]_t)[y - x] \| & \leq & 
\frac{H_{\nu} (1 - t)^{1 + \nu}}{1 + \nu} \|y - x\|^{1 + \nu}, \\
\\
\| A(x) - A([x, y]_t) - tDA([x, y]_t)[x - y] \| & \leq & 
\frac{H_{\nu} t^{1 + \nu}}{1 + \nu} \|y - x\|^{1 + \nu}.
\ea
\eeq
Multiplying them by $t$ and $1 - t$, correspondingly, and using triangle inequality, we obtain
$$
\ba{rcl}
\frac{1}{t (1 - t)}
\| tA(y) + (1 - t)A(x) - A([x, y]_t) \| & \overset{(\ref{TwoHolderBounds})}{\leq} &
\frac{H_{\nu}}{1 + \nu} \|y - x\|^{1 + \nu} \cdot \Bigl[ t^{\nu} + (1 - t)^{\nu}  \Bigr] \\
\\
& \leq &
\frac{2^{1 - \nu} H_{\nu}}{1 + \nu} \|y - x\|^{1 + \nu} 
\ea
$$
Hence, we get that $\Delta_A(x, y) \leq \frac{2^{1 - \nu} H_{\nu}}{1 + \nu} \|y - x\|^{1 + \nu}$,
which leads to~\eqref{HolderGCB}.
\qed

\newpage
\section{Curvature of High-Order Derivatives}
\label{SectionHigh}

From now on, the main working metric space for us will be a normed vector space (Example~\ref{ExampleVectorSpace}).

In what follows, we denote by $\E$ a real finite-dimensional vector space.
Then, $\E^{*}$ is the dual space, which is the space of linear functions on $\E$. The value of a function $s \in \E^{*}$ on a vector $x \in \E$ is denoted by $\la s, x \ra = s(x) \in \R$. 

We fix a positive definite self-adjoint operator $B: \E \to \E^{*}$, $B \succ 0$, and use the following pair of Euclidean norms:
$$
\ba{rcl}
\| x \| & := & \la  Bx, x \ra^{1/2}, \qquad \| s \|_* \;\; := \;\; \la s, B^{-1} s\ra^{1/2}, \qquad x \in \E, s \in \E^{*},
\ea
$$
that satisfy the standard Cauchy inequality: $| \la s, x \ra | \leq \| s\|_* \|x\|$.

Let $V(\cdot)$ be a mapping defined on an open domain $\dom V$ containing a closed convex set $Q \subseteq \E$, with values in $\E^{*}$:
$$
\ba{c}
V : \dom V \to \E^{*}, \qquad Q \subset \dom V \subseteq \E.
\ea
$$
We will assume that $V$ is differentiable and denote by $DV(x)[h] \in \E^*$ its directional derivative at point $x \in Q \subset \dom V$
along direction $h \in \E$. Note that we can treat it as a bi-linear form on $\E$:
$$
\ba{rcl}
\la DV(x)h_1, h_2 \ra & := & \la DV(x)[h_1], h_2 \ra \;\; \in \;\; \R, \qquad x \in Q \subseteq \E^*.
\ea
$$
However, this bi-linear form is not necessary symmetric.
Inductively, for any $k \geq 0$, we denote by $D^k V(x)[h_1, \ldots h_k] \in \E^*$ the $k$th directional derivative of $V$ along a fixed set of directions $h_1, \ldots h_k \in \E$. 

Note that $D^k V(x)$ is symmetric with respect to permuting the directions, as long as $V$ is sufficiently smooth.
When all directions are the same, $h_1 = \ldots = h_k = h \in \E$, we use a shorthand $D^k V(x)[h]^k := D^k V(x)[h, \ldots, h] \in \E^*$. 

In a natural way, we can treat $D^k V(x)$ as a multilinear form.
We also use notation $D^k V(x)[h]$ for the $(k-1)$th linear operator defined by applying $D^k V(x)$ to one fixed direction $h \in \E$.
For multilinear forms, we define the corresponding norm, as follows:
$$
\ba{rcl}
\| D^k V(x) \| & := & \max\limits_{\substack{h_1, \ldots h_{k} \in \E \\ \| h_i \| \leq 1 \; \forall i}}
\| D^k V(x)[h_1, \ldots, h_k] \|_*.
\ea
$$

Therefore, we endow the vector space of multilinear forms with the induced operator norm, and we can directly apply the results of the previous section to this metric space. Thus, for $k \geq 0$, the curvature of mapping $D^k V$, according to our definitions, is
$$
\ba{cl}
&\kappa_{D^k V}(r) \\[10pt]
& = 
\sup\limits_{\substack{x, y \in Q \\ 0 < t < 1}}
\Bigl\{
\frac{1}{t(1 - t)} \| t D^k V(y) + (1 - t) D^k V(x)
- D^k V( x + t (y - x) ) \| \; : \; \|y - x\| \leq r
\Bigr\}.
\ea
$$

In this paper, we are interested to study efficient algorithms for solving optimization and variational problems involving $V$, and describe their complexity via the exact GCB of different order:
\beq \label{CurvaturesK}
\ba{c}
\kappa_{V}(r), \; \kappa_{DV}(r), \; \kappa_{D^2V}(r), \; \ldots
\ea
\eeq
We will show that the curvature bounds~\eqref{CurvaturesK} serves as the main characteristic of the global complexity for the corresponding problem with $V$.
Note that we do not fix in advance a particular problem class for $V$, such as Lipschitz continuity of the derivative of a certain order. Our method does not need the problem class to be fixed as well.
We need to fix only
\textit{the order $p \geq 1$} of the derivatives that will be used in computations.
Therefore, our analysis allows for much more flexibility in terms of the local and global behavior of the operator $V$. At the same time, we demonstrate how to recover the correct complexity bounds in all known particular cases of the problem classes.

For $p \geq 0$, let us introduce the Taylor polynomial  of mapping $V$ at point $x$, as follows:
\beq \label{TaylorPolynomial}
\ba{rcl}
\mathcal{T}_{x, V}^p (y) & \Def &
\sum\limits_{k = 0}^p \frac{1}{k!} D^k V(x)[y - x]^k
\;\; \in \;\; \E^{\star}, \quad y \in Q \subseteq \E.
\ea
\eeq
For $p = 1$, it means that our method is \textit{first-order} in $V$,
and by Remark~\ref{RemarkLinearBound}, we immediately obtain the following global bound
on the linear approximation of $V$, using its GCB characteristic:
\beq \label{FirstOrderBound}
\ba{rcl}
\| V(y) - V(x) - DV(x)[y - x] \| & \leq & \kappa_{V}(\|y - x\|), \qquad x, y \in Q.
\ea
\eeq
For an arbitrary $p \geq 2$, we obtain global bounds on the Taylor approximation, using
the \textit{integral smoothing} of the GCB.
Define,
$$
\ba{rcl}
\hat{\sigma}_0(r) & \Def & \kappa_V(r),
\ea
$$
and, for $q \geq 1$, we denote
\beq \label{SigmaQDef}
\ba{rcl}
\hat{\sigma}_q(r) & \Def &
\frac{r^q}{(q - 1)!} \int\limits_0^1 (1 - \tau)^{q - 1} \kappa_{ D^q V } (\tau r) d\tau
\;\; = \;\;
\frac{1}{(q - 1)!} \int\limits_0^r (r - t)^{q - 1}
\kappa_{ D^q V } (t) dt.
\ea
\eeq
Note that $\hat{\sigma}_{q}(\cdot)$
are non-decreasing functions on its domain, and it holds, 
\beq \label{SigmaDer}
\ba{rcl}
\hat{\sigma}_1'(r) & = & \kappa_{DV}(r), \\
\\
\hat{\sigma}_{q}'(r) & = & 
\frac{1}{(q - 2)!} \int\limits_{0}^r (r - t)^{q - 2} \kappa_{D^q V}(t) dt
\;\; = \;\;
\frac{r^{q - 1}}{(q - 2)!}
\int\limits_0^1 (1 - \tau)^{q - 2} \kappa_{D^q V}(\tau r) d \tau,
\ea
\eeq
for $q  \geq 2$.
Since $\hat{\sigma}_{q}'(\cdot)$ is non-decreasing, we conclude that $\hat{\sigma}_{q}(\cdot), q \geq 1$ is \textit{convex}.

\begin{lemma} \label{LemmaVUp}
	For all $x, y \in Q$ and $p \geq 1$, we have
	\beq \label{VUp}
	\ba{rcl}
	\| V(y) - \mathcal{T}_{x, V}^{p}(y) \|_*
	& \leq & 
	\hat{\sigma}_{p - 1}(\|y - x\|).
	\ea
	\eeq
	Moreover, for $p \geq 2$, we have
	\beq \label{VDerUp}
	\ba{rcl}
	\| DV(y) - D \mathcal{T}_{x, V}^{p}(y) \| & \leq & 
	\hat{\sigma}_{p - 1}'( \| y - x \| ).
	\ea
	\eeq
\end{lemma}
\proof
The case $p = 1$ follows immediately from~\eqref{FirstOrderBound}.
Let us assume that $p \geq 2$.
Denote $r = \|y - x\|$.

Note that, by Taylor formula, for an arbitrary differentiable operator $A(\cdot)$ defined on $Q$ and for $q \geq 2$, we have
$$
\ba{rcl}
A(y) - \mathcal{T}_{x, A}^{q - 2}(y)
& = & 
\frac{1}{(q - 2)!}
\int\limits_0^1 (1 - \tau)^{q - 2} D^{q - 1} A(x + \tau (y - x))[y - x]^{q - 1} d\tau.
\ea
$$
We can extend this representation as follows:
\beq \label{TaylorExtended}
\ba{rcl}
& & \!\!\!\!\!\!\!\!\!\!\!\!\!\!\!\!\!\!\!\!\!\!
\| A(y) - \mathcal{T}_{x, A}^{q}(y) \|_* \\
\\
& \!\!\!\!\!\! = \!\!\!\!\!\! & 
\| A(y) - \mathcal{T}_{x, A}^{q - 2}(y)
- \frac{1}{(q - 1)!} D^{q - 1} A(x)[y - x]^{q - 1} 
- \frac{1}{q!} D^{q} A(x)[y - x]^{q} \|_* \\ 
\\
& \!\!\!\!\!\! = \!\!\!\!\!\! &
\| \frac{1}{(q - 2)!}
\int\limits_{0}^1 (1 - \tau)^{q - 2} \Bigl[ D^{q - 1} A(x + \tau(y - x)) 
\\
\\
& & \quad - D^{q - 1} A(x) - \tau D^{q} A(x)[y - x]
\Bigr] [y - x]^{q - 1} d\tau \|_* \\
\\
&  \!\!\!\!\!\! \overset{(\ref{eq-AppUp})}{\leq} \!\!\!\!\!\! &
\frac{r^{q - 1}}{(q - 2)!} \int\limits_{0}^{1}(1 - \tau)^{q - 2} \kappa_{D^{q - 1} A}(\tau r) d \tau.
\ea
\eeq
Therefore, to prove~\eqref{VUp}, it is enough to substitute $A(\cdot) := V(\cdot)$, $q := p$,
and to notice that the right hand side of~\eqref{TaylorExtended} is $\hat{\sigma}_{p - 1}(r)$.

Now, let us prove~\eqref{VDerUp}. The case $p = 2$ is immediate, as
$$
\ba{rcl}
\| DV(y) - D \mathcal{T}_{x, V}^2(y) \|
& = &
\| DV(y) - DV(x) - D^2V(x)[y - x] \| \\
\\
& \overset{(\ref{eq-AppUp})}{\leq} &
\kappa_{DV}(\|y - x\|)
\;\; \overset{(\ref{SigmaDer})}{=} \;\;
\hat{\sigma}_{1}'(\|y - x\|).
\ea
$$
The case $p \geq 3$ follows by substituting $A(\cdot) := DV(\cdot)$ and $q := p - 1$ into~\eqref{TaylorExtended}, which gives
$$
\ba{rcl}
\| DV(y) - D\mathcal{T}_{x, V}^p(y) \|
& = & 
\| DV(y) - \mathcal{T}_{x, DV}^{p - 1}(y) \| \\
\\
& \overset{(\ref{TaylorExtended})}{\leq} &
\frac{r^{p - 2}}{(p - 3)!} \int\limits_0^1 (1 - \tau)^{p - 3} \kappa_{D^{p - 1} V}(\tau r) d\tau
\;\; \overset{(\ref{SigmaDer})}{=} \;\;
\hat{\sigma}'_{p - 1}(r).
\ea
$$
Thus,~\eqref{VDerUp} is established for all $p \geq 2$.
\qed

To conclude this section, let us list the most important properties of $\hat{\sigma}_q(\cdot)$, $q \geq 0$, that will be
crucial for our analysis.

\begin{lemma} \label{LemmaHSig}
	Let $q \geq 1$.
	For any $r \in \Gamma$ and $\beta \in [0, 1]$, we have
	\beq \label{SigLow}
	\ba{rcl}
	\hat{\sigma}_q(\beta r) & \geq & \beta^{q + 2} \hat{\sigma}_q(r),
	\qquad
	\hat{\sigma}_q'(\beta r) \;\; \geq \;\; \beta^{q + 1} \hat{\sigma}_q'(r).
	\ea
	\eeq
	Moreover,
	\beq \label{Compat}
	\ba{rcl}
	\hat{\sigma}_q(r) & \leq & r \hat{\sigma}_q'(r) 
	\;\; \leq \;\; (q + 2) \hat{\sigma}_q(r).
	\ea
	\eeq
\end{lemma}
\proof
Inequalities~\eqref{SigLow} follow from Theorem~\ref{TheoremMain}.
Since $\hat{\sigma}_q(0) = 0$, the first inequality in~\eqref{Compat}
follows from convexity of $\hat{\sigma}_q(\cdot)$. Finally, for small $\epsilon > 0$, we have
$$
\ba{rcl}
\hat{\sigma}_{q}(r) - \hat{\sigma}_{q}(r - \epsilon)
& = &
\hat{\sigma}_q(r) - \hat{\sigma}_{q}\bigl( r \cdot (1 - \frac{\epsilon}{r}) \bigr)
\;\; \overset{(\ref{SigLow})}{\leq} \;\;
\hat{\sigma}_q(r) - \bigl(1 - \frac{\epsilon}{r} \bigr)^{q + 2} \hat{\sigma}_q(r).
\ea
$$
Dividing this inequality by $\epsilon$ and taking the limit $\epsilon \to 0$, we get
$\hat{\sigma}_q'(r) \leq \frac{q + 2}{r} \hat{\sigma}_q(r)$.
\qed

\begin{corollary} \label{CorollaryLq}
	For any $q \geq 0$, the functions
	$$
	\ba{rcl}
	L_q(r) & := & \frac{\hat{\sigma}_{q}(r)}{r^{q + 2}}
	\qquad \text{and} \qquad
	\bar{L}_{q + 1}(r) \;\; := \;\; \frac{\hat{\sigma}_{q + 1}'(r)}{r^{q + 2}}
	\ea
	$$
	are non-increasing on their domain.
\end{corollary}
\proof
Indeed, for any $\beta \in [0, 1]$, due to~\eqref{KappaGrowth} and \eqref{SigLow} we have
$$
\ba{rcl}
L_q(\beta r)
& = & 
\frac{\hat{\sigma}_{q}(\beta r)}{(\beta r)^{q + 2}}
\;\; \geq \;\;
\frac{\hat{\sigma}_q(r)}{r^{q + 2}}
\;\; = \;\;
L_q(r),
\ea
$$
which concludes the proof for $L_q(\cdot)$. The proof for $\bar{L}_{q + 1}(\cdot)$ is analogous.
\qed

\begin{corollary} \label{CorollaryUp2}
	Let $q \geq 1$. For any $r, s \in \Gamma$, it holds
	\beq \label{SigmaPrimeUp2}
	\ba{rcl}
	\hat{\sigma}_{q}'(r) & \leq & \hat{\sigma}_{q}'(s)
	+ \frac{\hat{\sigma}_{q}'(s)}{s^{q + 1}} r^{q + 1}.
	\ea
	\eeq
\end{corollary}
\proof
Indeed, since $\hat{\sigma}_{q}'(\cdot)$ is non-decreasing, for $r \leq s$ we have
$\hat{\sigma}_{q}'(r) \leq \hat{\sigma}_{q}'(s)$. Otherwise, for $r > s$ we conclude by Corollary~\ref{CorollaryLq}
that
$$
\ba{rcl}
\frac{\hat{\sigma}_q'(r)}{r^{q + 1}} & \leq & \frac{\hat{\sigma}_q'(s)}{s^{q + 1}}
\qquad \Leftrightarrow \qquad
\hat{\sigma}_q'(r) \;\; \leq \;\; 
\frac{\hat{\sigma}_q'(s)}{s^{q + 1}} r^{q + 1}.
\ea
$$
Combining these two cases gives~\eqref{SigmaPrimeUp2}.
\qed

Finally, let us consider the following example, which enables us to bound 
the growth of the smoothed GCB characteristic $\hat{\sigma}_q(\cdot)$,
under H\"older continuity of certain high-order derivatives.

\begin{example} \label{ExampleHolderHigh}
	Let $p \geq 1$ and assume that
	the $p$th derivative of $V$ is H\"older continuous of degree $\nu \in [0, 1]$, with a constant $H_{p, \nu} > 0$:
	\beq \label{HolderHighDer}
	\ba{rcl}
	\| D^p V(x) - D^p V(y) \| & \leq & H_{p, \nu} \|x - y\|^{\nu}, \qquad x, y \in Q.
	\ea
	\eeq
	Then,
	\beq \label{SigmaHolderHigh}
	\ba{rcl}
	\hat{\sigma}_{p - 1}(r) & \leq & \frac{2^{1 - \nu}}{(1 + \nu) \cdot \ldots \cdot (p + \nu)} H_{p, \nu} r^{p + \nu}.
	\ea
	\eeq
\end{example}
\proof
Indeed, by Example~\ref{ExampleHolder}, we have
\beq \label{KappaHolderBound}
\ba{rcl}
\kappa_{D^{p - 1}V}(r) & \overset{(\ref{HolderGCB})}{\leq} & 
\frac{2^{1 - \nu} H_{p, \nu}}{1 + \nu} r^{1 + \nu}.
\ea
\eeq
This immediately gives~\eqref{SigmaHolderHigh} for $p = 1$, as $\hat{\sigma}_0(r) \equiv \kappa_{V}(r)$. 
Consider $p \geq 2$. We have
$$
\ba{rcl}
\hat{\sigma}_{p - 1}(r) & \overset{(\ref{SigmaQDef})}{=} & 
\frac{r^{p - 1}}{(p - 2)!} \int\limits_{0}^{1} (1 - \tau)^{p - 2} \kappa_{D^{p - 1}V}(\tau r) d \tau \\
\\
& \overset{(\ref{KappaHolderBound})}{\leq} &
\frac{r^{p + \nu}}{(p - 2)!} \cdot \frac{2^{1 - \nu} H_{p, \nu}}{1 + \nu}
\cdot I(p, \nu),
\ea
$$
where $I(p, \nu) := \int\limits_0^1 (1 - \tau)^{p - 2} \tau^{1 + \nu} d\tau = \frac{(p - 2)!}{(2 + \nu) \cdot \ldots \cdot (p + \nu)}$,
which gives~\eqref{SigmaHolderHigh}.
\qed

\section{Universal Reduced-Operator Method}
\label{SectionReduced}

Recall that we work with a mapping $V: Q \to \E^{*}$, where $Q \subseteq \E$ is a closed convex set.
Our target problem is the 
\textit{Composite Variational Ineqiality}:
\beq \label{ProblemCVI}
\ba{rcl}
\text{Find} \quad x^{\star} \in Q  \quad \text{s.t.} \quad
\la V(x^{\star}), x - x^{\star} \ra + \psi(x) & \geq & \psi(x^{\star}), \qquad \forall x \in Q,
\ea
\eeq
where $\psi(\cdot)$ is a closed convex function with domain $Q = \dom \psi$.
We assume that the mapping $V$ is several times differentiable and it constitutes the main difficulty in solving the problem~\eqref{ProblemCVI}. At the same time, $\psi$ has a sufficiently simple structure (e.g. an indicator of a convex set $Q$), for which we can efficiently compute the projection operation, $\proj_Q(y) := \argmin_{x \in Q} \|y - x\|$.

In this section, we develop our main algorithm for solving composite variational inequality~\eqref{ProblemCVI}
using \textit{reduced-operator} steps.

Our main assumption on $V$ is that for some  $p \geq 1$, the global curvature bound $\kappa_{D^{p - 1}V}(\cdot)$ 
is well defined and sufficiently smooth.
Namely, we require the following condition.

\begin{assumption} \label{AssumptionMonotone}
	For some $p \geq 1$, the function $\hat{\sigma}_{p - 1}(\cdot)$ is well-defined and strictly monotone on its domain.
\end{assumption}

From previous section, 
we already know that $\hat{\sigma}_{p - 1}(\cdot)$ is always non-decreasing, and it is convex when $p \geq 2$.
Thus, Assumption~\ref{AssumptionMonotone} on strict monotonicity is only needed to exclude some pathological cases,
as e.g. of the zero curvature $\sigma_{p - 1}(\cdot) \equiv 0$. Hence, by our assumption, 
the inverse of the GCB, $\sigma_{p - 1}^{-1}(\cdot)$,
is well-defined, and we will use this function to describe the global complexity of our algorithm.

The most important example of problem~\eqref{ProblemCVI} is when the vector field $V(\cdot)$ is \textit{potential}, i.e. it is the gradient of a certain differentiable function: $V \equiv \nabla f$ for 
some $f: \dom f \to \R$. In this case, \eqref{ProblemCVI}~is a minimization problem: 
$$
\ba{c}
\min_x F(x), \quad \text{with} \quad F(x) := f(x) + \psi(x).
\ea
$$
The exact GCB of the objective, $\kappa_f(\cdot)$, was used in~\cite{nesterov2025universal}
for justifying the universal complexity bounds for gradient methods, as applied to minimizing the composite function $F$.
Another important class of variational inequalities~\eqref{ProblemCVI} are min-max problems (see~\cite{nesterov2023high}).

Let us fix $p \geq 1$ (the order of the method). Hence, the method of order $p = 1$ means that we use the following linear approximation
that includes the Jacobian of the mapping:
$$
\ba{rcl}
V(y) & \approx & V(x) + DV(x)[y - x].
\ea
$$
For optimization problems, $V(x) \equiv \nabla f(x)$, and $p \geq 1$ means that we involve at least the Hessian of the objective, $DV(x) \equiv \nabla^2 f(x)$.

For a fixed $x \in Q$, $\alpha > 0$ and $M > 0$, denote the point $x^+ = x^+_p(x, \alpha, M)$
that satisfies the following variational conditions,
\beq \label{StatCond}
\ba{c}
\begin{cases}
	&V(x)
	+ \sum\limits_{k = 1}^{p} \frac{1}{k!} D^k V(x)[x^+ - x]^k
	+ (\alpha + Mr^{p}) B(x^+ - x) \; \in \; - \partial \psi(x^+), \\
	&r \; = \; \|x^+ - x\|.
\end{cases}
\ea
\eeq
Note that for $p = 1$ this step can be efficiently implemented by the means of Convex Optimization,
and for $p > 1$ it provides us with a higher-order generalization.
In the next section, we will show that for an appropriate choice of regularization parameters $\alpha, M > 0$,
we can ensure that the nonlinear operator in~\eqref{StatCond} is \textit{monotone}
and thus the subproblem is globally solvable.

The inclusion~\eqref{StatCond} gives us a particular choice of the subgradient $\psi'(x^+) \in \partial \psi(x^+)$. Let us denote
\beq \label{ReducedGradient}
\ba{rcl}
V_{\psi}(x^+) & \!\!\!\! := \!\!\!\! & V(x^+) + \psi'(x^+) \\
\\
& \!\!\!\! \overset{(\ref{StatCond})}{=} \!\!\!\! &	
V(x^+)
- V(x) - \sum\limits_{k = 1}^{p} \frac{1}{k!} D^k V(x)[x^+ - x]^k - (\alpha + Mr^{p})B(x^+ - x),
\ea
\eeq
which we call the \textit{Reduced Operator}. 
Note that by our definition, the operator $V_{\psi}(\cdot)$ is defined only at the points which are solutions to subproblem~\eqref{StatCond}.

The next lemma is the main technical tool, which provides us with the progress of the operation~\eqref{StatCond}, enabling us to establish universal convergence rates for our method.

\begin{lemma} \label{LemmaReducedGrad}
	Let $p \geq 1$,  $\delta > 0$, and let $M > 0$ be sufficiently large for ensuring the following condition:
	\beq \label{MLarge}
	\ba{rcl}
	\hat{\sigma}_{p - 1}\Bigl( 
	\Bigl[  \frac{2\delta}{5M} \Bigr]^{\frac{1}{p + 1}}
	\Bigr) & \leq & \frac{\delta}{5}.
	\ea
	\eeq
	Define $\alpha > 0$ as follows:
	\beq \label{AlphaChoice}
	\ba{rcl}
	\alpha & = & 
	M^{\frac{1}{p + 1}} \bigl[ \frac{2}{5} \delta \bigr]^{\frac{p}{p + 1}}.
	\ea
	\eeq
	Let $x_+$ be a solution of the variational inequality (\ref{StatCond}). If $\| V_{\psi}(x^+) \|_* \geq \delta$, then we have  
	\beq \label{RLowerBound}
	\ba{rcl}
	r & \geq & \Bigl[  \frac{2\delta}{5M} \Bigr]^{\frac{1}{p + 1}}
	\ea
	\eeq
	and
	\beq \label{MainProgress}
	\ba{rcl}
	\la V_{\psi}(x^+), x - x^+ \ra & \geq & 
	c_p \cdot 
	\Bigl(  \frac{ 1 }{M}   \| V_{\psi}(x^+) \|_*^{p + 2} \Bigr)^{\frac{1}{p + 1}},
	\ea
	\eeq
	where $c_p \Def  {1 \over 4}\left( {3 \over 2} \right)^{p \over 2(p+1)} \left[ \left({p+2 \over p} \right)^{p \over 2(p+1)} + \left({p \over p +2} \right)^{p +2 \over 2(p+1)}\right] \geq {1 \over 3}$ is a numerical constant.
\end{lemma}
\proof
Indeed, by definition of the step~\eqref{StatCond} and by Taylor approximation~\eqref{VUp}, we have 
\beq \label{EuclideanBound}
\ba{rcl}
\!\!\!\!\!\!
\| V_{\psi}(x^+) + (\alpha + Mr^{p}) B(x^+ - x) \|_*
& \! \overset{(\ref{StatCond})}{=} \! &
\| V(x^+) - \mathcal{T}_{x, V}^{p}(x^+) \|_*
\; \overset{(\ref{VUp})}{\leq} \;
\hat{\sigma}_{p - 1}( r ). 
\ea
\eeq
Therefore, by triangle inequality,
\beq \label{Triangle}
\ba{rcl}
\delta & \leq & \| V_{\psi}(x^+) \|_*
\;\; \leq \;\; 
\alpha r + Mr^{p + 1} + \hat{\sigma}_{p - 1}(r).
\ea
\eeq
Now, assume that
\beq \label{ContrAssumption}
\ba{rcl}
r & < & \Bigl[ \frac{2 \delta}{5M} \Bigr]^{\frac{1}{p + 1}}
\qquad 
\Leftrightarrow
\qquad
M r^{p + 1} \;\; < \;\; \frac{2\delta}{5}.
\ea
\eeq
It also leads to
$$
\ba{rcl}
\alpha r & = & 
M^{\frac{1}{p + 1}}
\bigl[  \frac{2}{5} \delta \bigr]^{\frac{p}{p + 1}}
r
\;\; \overset{(\ref{ContrAssumption})}{<} \;\;
\frac{2 \delta}{5}.
\ea
$$
Then, according to~\eqref{Triangle}, 
it holds $\hat{\sigma}_{p - 1}(r) > \frac{\delta}{5}$,
which combined with~\eqref{MLarge} leads to
$$
\ba{rcl}
\hat{\sigma}_{p - 1}\Bigl( \Bigl[  \frac{2\delta}{5M} \Bigr]^{\frac{1}{p + 1}} \Bigr)
& < & \hat{\sigma}_{p - 1}(r) 
\qquad \Rightarrow \qquad
\Bigl[  \frac{2\delta}{5M} \Bigr]^{\frac{1}{p + 1}}
\;\; < \;\; r,
\ea
$$
due to monotonicity of $\hat{\sigma}_{p - 1}(\cdot)$,
that contradicts~\eqref{ContrAssumption}. Therefore, \eqref{RLowerBound} is proved.

Denote $\beta := \alpha + M r^p$.
In view of Corollary~\ref{CorollaryLq}, function 
$L(r) = \frac{ \hat{\sigma}_{p - 1}(r) }{r^{p + 1}}$ is non-increasing in~$r$. Hence,
$$
\ba{rcl}
\frac{ \hat{\sigma}_{p - 1}(r) }{r^{p + 1}} & \overset{(\ref{RLowerBound})}{\leq} &
\frac{ \hat{\sigma}_{p - 1}\bigl(  \bigl[  \frac{2\delta}{5M} \bigr]^{1 / (p + 1)}   \bigr) }{   \frac{2\delta}{5M}  }
\;\; \overset{(\ref{MLarge})}{\leq} \;\;
\frac{M}{2},
\ea
$$
and therefore
$\hat{\sigma}_{p - 1}(r) \leq \half 
Mr^{p+1} \leq \half r \beta$.
Taking the square of both sides of~\eqref{EuclideanBound} , we get
\beq \label{InnerProdProgress}
\ba{rcl}
\la V_{\psi}(x^+), x - x^+ \ra & \overset{(\ref{EuclideanBound})}{\geq} & 
\frac{1}{2 \beta } \| V_{\psi}(x^+) \|_*^2
+ \frac{ \beta r^2}{2}
- \frac{\hat{\sigma}_{p - 1}(r)^2}{2 \beta} \\ 
\\

& \geq & 
\frac{1}{2 \beta} \| V_{\psi}(x^+) \|_*^2
+ \frac{3 }{8} \beta r^2.
\ea
\eeq
Since $\alpha \overset{(\ref{RLowerBound})}{\leq} M r^p$, we obtain
\beq \label{InnerProdProgress2}
\ba{rcl}
\la V_{\psi}(x^+), x - x^+ \ra
& \geq & 
\frac{1}{2( \alpha + Mr^p )} \| V_{\psi}(x^+) \|_*^2
+ \frac{3}{8} (\alpha + Mr^p) r^2\\
\\
& \geq &
\frac{1}{4Mr^p} \| V_{\psi}(x^+) \|_*^2 
+ \frac{3}{8}Mr^{p + 2}.
\ea
\eeq
The direct minimization of the right hand side of~\eqref{InnerProdProgress2} in $r > 0$ gives us (\ref{MainProgress}). It remains to note that
$$
\ba{rcl}
c_p & \geq & \hat c_p \Def \half \left( {3 \over 2} \right)^{p \over 2(p+1)} \left({p \over p +2} \right)^{1 \over 2(p+1)} \; \geq \; \hat c_1 \; = \; \left( \half \right)^{5/4} \; \geq  \; {1 \over 3}. \hspace{5ex} \rule{1.0ex}{1.0ex}
\ea
$$

We are ready to present our universal algorithm for solving the variational problem~\eqref{ProblemCVI}.
It can use any of the following \textit{convergence conditions}, which are
classical for solving optimization and variational problems:
\begin{itemize}
	\item A small norm of the Reduced Operator: $\| V_{\psi}( x ) \|_* \leq \delta$, for a certain $\delta > 0$;
	
	\item A small value of an Accuracy Certificate: $\Delta_k \leq \varepsilon$, for a certain $\epsilon > 0$. 
\end{itemize}
Our accuracy certificates $\Delta_k$ will be constructed algorithmically during iterations $k \geq 1$, 
and will provide us with a computable upper bound to the following \textit{merit function}, 
$$
\ba{rcl}
\mu(\bar{x}) & \Def & 
\psi(\bar{x}) + \max\limits_{x \in Q} \Bigl[ 
\la V(x), \bar{x} - x\ra - \psi(x)
\Bigr],
\ea
$$
which is a standard measure of accuracy for problem~\eqref{ProblemCVI} ,
when the operator $V$ is monotone (see~\cite{nesterov2023high}).

We will ensure that $\Delta_k \geq \mu(\bar{x}_k)$ for a certain sequence of points $\{ \bar{x}_k \}_{k \geq 1}$ and
accuracy certificates $\{ \Delta_k \}_{k \geq 1}$, generated within our algorithm. A technical assumption for
computing $\Delta_k$ is that the domain $Q$ of our problem is bounded.
This is not restrictive, as we can always introduce an auxiliary constraint to our problem of the form $\| x \| \leq D$, 
ensuring that $D$ is sufficiently large, so the constraint will be satisfied for the solution to our initial problem.

Let us describe an explicit structure of the certificates $\Delta_k$.
For a sequence of test points $\{ x_k \}_{k \geq 1}$ generated by the method and corresponding sequence of evaluated reduced operators~$\{ V_{\psi}(x_k) \}_{k \geq 1}$, we denote
$$
\ba{rcl}
\Delta_k & := & 
\frac{1}{A_k} \max\limits_{x \in Q} \sum\limits_{i = 1}^k a_i \la V_{\psi}(x_i), x_i - x \ra,
\ea
$$
where $\{ a_i \}_{i \geq 1}$ are positive weights and $A_k := \sum_{i = 1}^k a_i$.
Now, assuming that operator $V(\cdot)$ is monotone,
\beq \label{VMonotone}
\ba{rcl}
\la V(y) - V(x), y - x \ra & \geq & 0, \qquad x, y \in Q,
\ea
\eeq
and taking into account convexity of $\psi$, we have
$$
\ba{rcl}
\Delta_k & \overset{(\ref{ReducedGradient})}{\geq} &
\frac{1}{A_k} \max\limits_{x \in Q} 
\sum\limits_{i = 1}^k a_i \Bigl[ \la V(x_i), x_i - x \ra + \psi(x_i) - \psi(x) \Bigr] \\[10pt]
& \overset{(\ref{VMonotone})}{\geq} &
\frac{1}{A_k} \max\limits_{x \in Q} 
\sum\limits_{i = 1}^k a_i \Bigl[ \la V(x), x_i - x \ra + \psi(x_i) - \psi(x) \Bigr] \\[10pt]
& \geq & 
\max\limits_{x \in Q} \Bigl[ \la V(x), \bar{x}_k - x \ra + \psi(\bar{x}_k) - \psi(x) \Bigr]
\;\; = \;\; \mu(\bar{x}_k),
\ea
$$
where we defined the average point as follows:
$$
\ba{rcl}
\bar{x}_k & := & \frac{1}{A_k} \sum\limits_{i = 1}^k x_i.
\ea
$$
Note that for the potential case (i.e., $V \equiv \nabla f$ for a convex function $f$), we have
$$
\ba{rcl}
\mu(\bar{x}_k) & \leq & F(\bar{x}_k) - F^{\star} \;\; \leq \;\; \Delta_k.
\ea
$$
We also note that the guarantee of the small reduced operator norm is stronger, as it implies a small value of the merit function at \textit{the last point}. Indeed, 
by the Cauchy inequality:
$$
\ba{rcl}
\| V_{\psi}(x_k) \|_*
& \!\!\! \geq \!\!\! & \frac{1}{D} \max\limits_{x \in Q} \la V_{\psi}(x_k), x_k - x \ra 
\overset{(\ref{ReducedGradient})}{\geq}
\frac{1}{D} \max\limits_{x \in Q} 
\Bigl[ \la V(x_k), x_k - x \ra + \psi(x_k) - \psi(x) \Bigr] \\[5pt]
& \!\!\! \overset{(\ref{VMonotone})}{\geq} \!\!\! & 
\frac{1}{D} \max\limits_{x \in Q} 
\Bigl[ \la V(x), x_k - x \ra + \psi(x_k) - \psi(x) \Bigr] 
\;\; = \;\;
\frac{1}{D} \mu(x_k).
\ea
$$
We use the accuracy in terms of the merit function as the main complexity measure of our method. 
However, the bound in terms of the reduced operator norm provides an additional stopping condition.

\beq\label{AlgorithmMain}
\ba{|c|}
\hline\\
\quad \mbox{\bf Universal Reduced-Operator Method, $p \geq 1$.} \quad\\[10pt]
\hline\\
\ba{l}
\mbox{{\bf Initialization.} Choose $x_0 \in Q$, $M_0 > 0$, and $\delta, \varepsilon > 0$. Set $v_0 = x_0$, $A_0 = 0$.} \\[7pt]
\mbox{\bf For $k \geq 0$ iterate:} \\[7pt]
\mbox{\textbf{1.} Find the smallest integer $i_k \geq 0$ s.t. for} \\[7pt]
\mbox{\qquad $M^+ := 2^{i_k}M_k$  \quad and 
	\quad $x^+ := x_p^+(v_k, (M^+)^{\frac{1}{p + 1}} \bigl[ \frac{2}{5} \delta \bigr]^{\frac{p}{p + 1}}, M^+)$  it holds:}\\[7pt]
\mbox{\qquad $\la V_{\psi}(x^+), v_k - x^+ \ra \geq  c_p \cdot \Bigl[ \frac{\| V_{\psi}(x^+) \|_*^{p + 2} }{M^+}  \Bigr]^{1 / (p + 1) }$
	\text{or}  $\| V_{\psi}(x^+) \|_* \leq \delta$.} \\[7pt]
\mbox{\textbf{2.} Set $x_{k + 1} = x^+$ and $M_{k + 1} = \frac{1}{2}M^+$.} \\[7pt]
\mbox{\textbf{3.} If $\| V_{\psi}(x_{k + 1}) \|_* \leq \delta$ then \textbf{stop} and \textbf{return} $x_{k + 1}$. } \\[9pt]
\mbox{\textbf{4.} Set  $a_{k + 1} = \frac{\la V_{\psi}(x_{k + 1}), v_k - x_{k + 1}\ra }{\| V_{\psi}(x_{k + 1}) \|_*^2}$,
	$A_{k + 1} = A_k + a_{k + 1}$,
	and compute } \\[7pt]
\mbox{\qquad $\Delta_{k + 1} := \frac{1}{A_{k + 1}} \max\limits_{x \in Q} \sum\limits_{i = 1}^{k + 1} a_i \la V_{\psi}(x_i), x_i - x \ra$.} \\[7pt]
\mbox{\textbf{5.} If $\Delta_{k + 1} \leq \varepsilon$ 
	then \textbf{stop} and \textbf{return}  $\bar{x}_{k + 1} = \frac{1}{A_{k + 1}} \sum\limits_{i = 1}^{k + 1} a_i x_i $.} \\[7pt]
\mbox{\textbf{6.} Perform the reduced-operator step: } \\[7pt]
\qquad
\ba{rcl}
v_{k + 1}  & = &  \proj_Q\bigl( 
v_k - a_{k + 1} B^{-1} V_{\psi}(x_{k + 1})\bigr).
\ea
\\
\ea
\\
\\
\hline
\ea
\eeq

We assume that there exists a point $x^{\star} \in Q$ (weak solution) satisfying the following condition:
\beq \label{WeakSolution}
\ba{rcl}
\la V(x), x - x^{\star} \ra +  \psi(x)  -  \psi(x^{\star})& \geq & 0, \qquad \forall x \in Q.
\ea
\eeq
Condition~\eqref{WeakSolution} is satisfied, e.g., when the operator $V$ is monotone~\eqref{VMonotone}
and $x^{\star}$ is a solution to~\eqref{ProblemCVI}. We denote $R_0 \Def \|x_0 - x^{\star}\| \leq D \Def \diam Q$.

\begin{theorem} \label{TheoremReducedGrad}
	Let the Assumption~\ref{AssumptionMonotone} hold and relation
	\eqref{WeakSolution} be satisfied for some $x^{\star} \in Q$.
	Assume that algorithm~\eqref{AlgorithmMain} does not stop for the first $K \geq 1$ iterations.
	Then,
	\beq \label{ComplexityBoundReducedGrad}
	\ba{rcl}
	K & \leq &
	\bigl( \frac{4}{5} \bigr)^{2 \over p + 1}
	\Bigl[ 
	\,
	\frac{R_0}{c_p}  \max\Bigl\{ 
	\frac{1}{ \hat{\sigma}_{p - 1}^{-1}( \delta / 5 ) },
	\Bigl[ \frac{5M_0}{2 \delta} \Bigr]^{1 \over p + 1}
	\Bigr\}
	\,
	\Bigr]^2
	\cdot
	\ea
	\eeq
	Assume additionally that the operator $V$ is monotone~\eqref{VMonotone}, and set $\delta = \frac{\varepsilon}{D}$.
	Then, 
	\beq \label{ComplexityAccuracyCertificate}
	\ba{rcl}
	K & \leq & \left(4 \over 5\right)^{2 \over p+2} \left[ {D \over c_p} \max \left\{ {1 \over \hat \sigma^{-1}_{p-1}(\varepsilon /(5D))}, \left[ 5 M_0 D \over 2 \epsilon \right]^{1 \over p+1} \right\} \right]^{2 (p+1) \over p+2}.
	\ea
	\eeq
\end{theorem}
\proof
Note that by Lemma~\ref{LemmaReducedGrad}, each iteration of our algorithm is well defined.
Consider iteration $0 \leq k \leq K - 1$.
By the definition of the projection operation, we have
$$
\ba{rcl}
v_{k + 1} & \! = \! & 
\argmin\limits_{x \in Q} \| v_k - a_{k + 1} B^{-1} V_{\psi}(x_{k + 1}) - x \|^2 \\[10pt]
& \! = \! & 
\argmin\limits_{x \in Q} \Bigl\{
\frac{1}{2}\|v_k - x\|^2 + a_{k + 1} \la V_{\psi}(x_{k + 1}), x  \ra 
\Bigr\}.
\ea
$$
Hence, due to strong convexity of this optimization subproblem, we have, for any $x \in Q$:
$$
\ba{cl}
& \frac{1}{2} \|v_k - x\|^2 + a_{k + 1} \la V_{\psi}(x_{k + 1}), x - x_{k + 1} \ra \\
\\
& \geq \;
\frac{1}{2} \|v_{k + 1} - x\|^2 + \frac{1}{2}\| v_{k} - v_{k + 1} \|^2 + a_{k + 1} \la V_{\psi}(x_{k + 1}), v_{k + 1} - x_{k + 1} \ra \\
\\
& = \;
\frac{1}{2} \|v_{k + 1} - x\|^2 + \frac{1}{2} \| v_k - v_{k + 1} \|^2 - a_{k + 1} \la V_{\psi}(x_{k + 1}), v_k - v_{k + 1} \ra 
+ a_{k + 1}^2 \| V_{\psi}(x_{k + 1}) \|^2_*,
\ea
$$
where we used the definition of $a_{k + 1}$ in the last equation.
Therefore, we get
$$
\ba{rcl}
\frac{1}{2} \|v_k - x\|^2 & \geq & \frac{1}{2} \|v_{k + 1} - x\|^2
+ a_{k + 1} \la V_{\psi}(x_{k + 1}), x_{k + 1} - x \ra + a_{k + 1}^2 \| V_{\psi}(x_{k + 1}) \|^2_* \\
\\
& &  \quad + \; \min\limits_{\tau > 0} \Bigl\{  \frac{1}{2}\tau^2 - a_{k + 1} \| V_{\psi}(x_{k + 1})\|_* \tau  \Bigr\} \\
\\
& = & 
\frac{1}{2} \|v_{k + 1} - x\|^2
+ a_{k + 1} \la V_{\psi}(x_{k + 1}), x_{k + 1} - x \ra + \frac{1}{2} a_{k + 1}^2 \| V_{\psi}(x_{k + 1}) \|^2_*.
\ea
$$
Telescoping this bound for the first $K$ iterations, we obtain
\beq \label{ReducedGradGuarantee}
\ba{rcl}
\frac{1}{2}\| x_0 - x \|^2 
& \geq & 
\frac{1}{2}\| v_K - x\|^2 
+ \frac{1}{2} \sum\limits_{i = 1}^K a_i^2 \| V_{\psi}(x_i) \|_*^2
+ \sum\limits_{i = 1}^K a_i \la V_{\psi}(x_i), x_i - x \ra.
\ea
\eeq

By the definition of our coefficients, we have the relationship:
\beq \label{MkRel}
\ba{rcl}
M_{k + 1} & = & 2^{i_k - 1} M_k.
\ea
\eeq
In case $i_k \geq 1$, we conclude that, due to Lemma~\ref{LemmaReducedGrad}, $M_{k + 1}$ satisfies inequality
\beq \label{MkNextBound}
\ba{rcl}
\hat{\sigma}_{p - 1}\Bigl( \Bigl[  \frac{2 \delta}{5M_{k + 1}} \Bigr]^{1 / (p + 1)}  \Bigr) & \geq & \frac{\delta}{5},
\ea
\eeq 
and, otherwise, $M_{k + 1} < M_k$. Therefore, we establish the following bound for all coefficients
\beq \label{MkBound}
\ba{rcl}
M_{k} & \leq & \mathcal{M}
\;\; \Def \;\;
\max\Bigl\{  \frac{2\delta}{ 5[ \hat{\sigma}_{p - 1}^{-1}(\delta / 5) ]^{p + 1} }, M_0  \Bigr\},
\qquad 0 \leq k \leq K.
\ea
\eeq

Note that by the definition of $a_{k + 1}$ and the condition on the progress in Step~1, we have
\beq \label{ReducedGradientakBound}
\ba{rcl}
a_{k + 1} & \geq & 
\frac{c_p}{(2M_{k + 1})^{1/(p + 1)}} \cdot \Bigl[ \frac{1}{\| V_{\psi}(x_{k + 1})\|_*} \Bigr]^{\frac{p}{p + 1}}
\;\; \overset{(\ref{MkBound})}{\geq} \;\;
\frac{c_p}{(2 \mathcal{M})^{1/(p + 1)}} \cdot \Bigl[ \frac{1}{\| V_{\psi}(x_{k + 1})\|_*} \Bigr]^{\frac{p}{p + 1}}.
\ea
\eeq
Therefore, substituting $x := x^{\star}$ into~\eqref{ReducedGradGuarantee}, and using that
$$
\ba{rcl}
\la V_{\psi}(x_i), x_i - x^{\star} \ra
& \overset{(\ref{ReducedGradient})}{\geq} &
\la V(x_i), x_i - x^{\star} \ra + \psi(x_i) - \psi(x^{\star})
\;\; \overset{(\ref{WeakSolution})}{\geq} \;\; 0, \qquad 1 \leq i \leq K,
\ea
$$
we obtain
$$
\ba{rcl}
R_0^2 & \geq &
\sum\limits_{i = 1}^K a_i^2 \| V_{\psi}(x_i) \|_*^2
\;\; \overset{(\ref{ReducedGradientakBound})}{\geq} \;\;
\frac{c_p^2}{(2 \mathcal{M})^{2 / (p + 1)}}
\sum\limits_{i = 1}^K \| V_{\psi}(x_i) \|_*^{2 / (p + 1)} \\[10pt]
& \geq &
\frac{c_p^2}{2^{2 / (p + 1)}} \cdot \Bigl[  \frac{\delta}{\mathcal{M}} \Bigr]^{2 / (p + 1)} \cdot K,
\ea
$$
where we used the stopping condition of our algorithm in the last inequality.
Using the definition of $\mathcal{M}$ immediately gives the required bound~\eqref{ComplexityBoundReducedGrad}.

To prove~\eqref{ComplexityAccuracyCertificate}, we notice that maximizing the left- and the right-hand side of~\eqref{ReducedGradGuarantee} 
with respect to $x \in Q$, we get:
\beq \label{DeltaKBound}
\ba{rcl}
D^2 & \geq & \sum\limits_{i = 1}^K a_i^2 \| V_{\psi}(x_i) \|_*^2 + 2 A_K \Delta_K  \;\; \overset{(\ref{ReducedGradientakBound})}{\geq} \;\;
\sum\limits_{i = 1}^K \left[ \zeta_p
a_i^{- 2 / p}
+2 a_i \Delta_K \right],
\ea
\eeq
where $\zeta_p \Def \frac{c_p^{2(p + 1) / p }}{ (2\mathcal{M})^{2 / p} }$.
Note that $\min\limits_{\tau>0} \left[ \zeta_p
\tau^{- 2 / p}
+2 \tau \Delta_K \right] =  {p+2 \over p} \zeta_p^{p \over p+2} \left (p \Delta_K\right)^{2 \over p+2}$. Assuming that $\Delta_K \geq \varepsilon$, we get the following bound:
$$
\ba{rcl}
K & \leq & {1 \over p+2} p ^{p \over p+2}  \frac{D^2}{\varepsilon^{\frac{2}{p + 2}} \zeta_p^{\frac{p}{p + 2}}} 
\;\; = \;\;
{1 \over p+2} p ^{p \over p+2}  \Bigl[ 
\frac{D}{c_p} \Bigl(  \frac{2 \mathcal{M} D}{\varepsilon} \Big)^{\frac{1}{p + 1}}
\Bigr]^{\frac{2(p + 1)}{p + 2}}
\\
\\
& = & {1 \over p+2} p ^{p \over p+2} \left(4 \over 5\right)^{2 \over p+2} \left[ {D \over c_p} \max \left\{ {1 \over \hat \sigma^{-1}_{p-1}(\varepsilon /(5D))}, \left[ 5 M_0 D \over 2 \varepsilon \right]^{1 \over p+1} \right\} \right]^{2 (p+1) \over p+2}.
\ea
$$
Since ${1 \over p+2} p ^{p \over p+2} \leq 1$, we come to the relation (\ref{ComplexityAccuracyCertificate}).
\qed

\begin{remark}
	Let the initial $M_0$ be sufficiently small, $M_0 \leq \frac{2 \delta}{5[ \hat{\sigma}_{p - 1}^{-1}(\delta / 5) ]^{p+1}}$, then the universal complexity of the method is as follows:
	$$
	\ba{c}
	O\Bigl( \Bigl[ \frac{R_0}{\hat{\sigma}_{p - 1}^{-1}(\delta / 5)} \Bigr]^2 \Bigr)
	\qquad \text{and} \qquad
	O\Bigl( 
	\Bigl[ 
	\frac{D}{\hat{\sigma}_{p - 1}^{-1}( \varepsilon / (5D) )}
	\Bigr]^{\frac{2(p + 1)}{p + 2}}
	\Bigr),
	\ea
	$$
	correspondingly, for reaching $\delta$-accuracy in terms of the reduced operator norm, and for 
	reaching $\varepsilon$-accuracy in terms of the merit function.
\end{remark}

\BR
For practical computations, it may be reasonable to choose $M_0$ as follows:
\beq\label{eq-M_0}
\ba{rcl}
M_0  & = &  {2 \varepsilon \over  5} \left( {1 \over c_p} \right)^{p+1}   D^p.
\ea
\eeq
Then, if method (\ref{AlgorithmMain}) does not stop after one iteration, then
\beq \label{eq-CompDelta}
\ba{rcl}
1 \; \leq \; K & \leq & \left(4 \over 5\right)^{2 \over p+2} \left[ {D \over c_p \, \hat \sigma^{-1}_{p-1}(\varepsilon /(5D))} \right]^{2 (p+1) \over p+2}.
\ea
\eeq
\ER

\begin{remark}
	To compute the accuracy certificate $\Delta_k$, we need
	an estimate $D \geq R_0$ on the diameter of the set $Q$.
	This is the only parameter of our method. It does not depend on the smoothness structure
	of the operator $V$. Thus, the algorithm is universal with respect to $\hat{\sigma}_{p - 1}(\cdot)$,
	which is the main characteristic of complexity in~\eqref{ComplexitiesSimplified}.
\end{remark}

\begin{corollary} \label{CorollaryRateHolder}
	For $p \geq 1$, assume that $D^p V$ is H\"older continuous of degree $\nu \in [0, 1]$
	with a constant $H_{p, \nu} > 0$. Then (see Example~\ref{ExampleHolderHigh}),
	$$
	\ba{rcl}
	\hat{\sigma}_{p - 1}(r) & \leq & c_{p, \nu} H_{p, \nu} r^{p + \nu}
	\qquad \Rightarrow \qquad
	\hat{\sigma}_{p - 1}^{-1}(y) \;\; \geq \;\; 
	\bigl[ 
	\frac{y}{c_{p, \nu} H_{p, \nu}}
	\bigr]^{\frac{1}{p + \nu}}
	\ea
	$$
	where $c_{p, \nu} := \frac{2^{1 - \nu}}{(1 + \nu) \cdot \ldots \cdot (p + \nu)}$ is a numerical constant.
	Therefore, we get the following pair of complexities for algorithm~\eqref{AlgorithmMain}:
	\beq \label{ComplexityPair}
	\ba{rcl}
	O\Bigl(  R_0^2 \cdot \Bigl[  \frac{H_{p, \nu}}{\delta} \Bigr]^{\frac{2}{p + \nu}}  \Bigr)
	\qquad & \text{and} & \qquad
	O\Bigl(  
	\Bigl[ 
	\frac{H_{p, \nu} D^{p + 1 + \nu}}{\varepsilon}
	\Bigr]^{\frac{2(p + 1)}{(p + 2)(p + \nu)}}
	\Bigr)
	\ea
	\eeq
	oracle calls, to get
	a point $x$ with the corresponding guarantee $\| V_{\psi}(x) \|_{*} \leq \delta$ and $\mu(x) \leq \varepsilon$.
	These bounds recover the bounds on the corresponding classes from~\cite{nesterov2023high}.
	At the same time, note that the method does not require 
	the knowledge of the actual structure of $\hat{\sigma}_{p - 1}(\cdot)$, automatically adjusting to the best problem class.
\end{corollary}

\begin{remark}
	For example, for $p = 1$ (first-order version of our algorithm) and considering
	the case $\nu = 1$ from Corollary~\ref{CorollaryRateHolder}, implying that the Jacobian $DV$ is Lipschitz with constant $L := H_{1, 1} > 0$, we obtain
	the pair of complexity bounds
	$$
	\ba{rcl}
	O \Bigl( \frac{L R_0^2}{\delta} \Bigr)
	\qquad & \text{and} & \qquad
	O \Bigl(  \frac{L^{2/3} D^2}{\varepsilon^{2/3}}  \Bigr).
	\ea
	$$
	We see that the complexity in terms of $\varepsilon$-accuracy on the merit function
	is better than that of $\delta$-accuracy on the reduced operator norm.
\end{remark}

To conclude this section, we note that in the potential case ($V \equiv \nabla f$ for a convex function $f$) with $p=1$ the performance guarantees for our method are slightly worse than that of the super-universal Newton method~\cite{doikov2024super}. However, our scheme is applicable to a much wider class of operators and it works for all $p \geq 1$. See also~\cite{grapiglia2022tensor} for the analysis of high-order tensor methods in the potential case.

\section{Monotonicity of the Subproblem}
\label{SectionMonotone}

In this section, we study whether the operator from our subproblem~\eqref{StatCond} is
\textit{monotone}, and thus it admits an efficient implementation. 
Clearly, it is monotone for $p = 1$, and hence we consider the nontrivial case $p \geq 2$.

We prove the following lemma.

\begin{lemma} \label{LemmaMonotone}
	Let operator $V$ be monotone and $p \geq 2$.
	Assume that $\delta > 0$ is fixed.
	Let us choose  sufficiently large $M > 0$,
	\beq \label{MBigMonotone}
	\ba{rcl}
	\hat{\sigma}_{p - 1} \Bigl( 
	\Bigl[ \frac{2 \delta}{5 M } \Bigr]^{\frac{1}{p + 1}}
	\Bigr) & \leq & \frac{2\delta}{5 (p + 1)},
	\ea
	\eeq
	and let us fix $\alpha > 0$ as previously:
	\beq \label{AlphaChoiceMonotone}
	\ba{rcl}
	\alpha & = & M^{\frac{1}{p + 1}}\bigl[ \frac{2}{5} \delta \bigr]^{\frac{p}{p + 1}}.
	\ea
	\eeq
	Then, the operator of subproblem~\eqref{StatCond},
	$$
	\ba{rcl}
	\Phi(y) & = & V(x) + \sum\limits_{k = 1}^p \frac{1}{k!} D^k V(x)[y - x]^k
	+ (\alpha + M \|y - x\|^p) B(y - x),
	\ea
	$$
	is monotone.
\end{lemma}
\proof
Note that
$$
\ba{rcl}
\Phi(y) & \equiv & \mathcal{T}_{x, V}^p(y)
+ \nabla d(y),
\ea
$$
where $d(\cdot)$ is the following the prox function,
$$
\ba{rcl}
d(y) & := & \frac{\alpha}{2}\|y - x\|^2 + \frac{M}{p + 2}\|y - x\|^{p + 2}.
\ea
$$
We have
$$
\ba{rcl}
\nabla d(y) & = & \alpha B(y - x) + M \|y - x\|^{p} B(y - x),
\ea
$$
and
\beq \label{ProxHess}
\ba{rcl}
\nabla^2 d(y)  & = & 
\alpha B +
M \|y - x\|^p B + 
p M \|y - x\|^{p - 2} B(y - x) [B(y - x)]^* \\
\\
& \succeq &
(\alpha + M \|y - x\|^p) B,
\ea
\eeq
where the operator $[a]^*: \E \to \R$ with $a \in \E^*$ works in accordance to the rule $[a]^* x = \la a, x \ra$, $x \in \E$.
Therefore,
$$
\ba{rcl}
\la D\Phi(y) h, h \ra
& = & 
\la DV(y)h, h \ra
+ \la ( D \mathcal{T}_{x, V}^p(y) - DV(y)) h, h \ra
+ \la \nabla^2 d(y) h, h \ra \\
\\
& \overset{(\ref{VDerUp}), (\ref{ProxHess})}{\geq} &
\la DV(y)h, h \ra
+ \Bigl(\alpha + M\|y - x\|^p - \hat{\sigma}_{p - 1}'(\|y - x\|) \Bigr) \| h \|^2 \\
\\
& \geq &
\Bigl(\alpha + M\|y - x\|^p - \hat{\sigma}_{p - 1}'(\|y - x\|) \Bigr) \| h \|^2,
\ea
$$
where the last inequality holds due to monotonicity of $V$.

To complete the proof it is enough to establish the following bound,
\beq \label{SufficientMonotone}
\ba{rcl}
\hat{\sigma}_{p - 1}'(\|y - x\|) & \leq & \alpha + M\|y - x\|^p.
\ea
\eeq
By Corollary~\ref{CorollaryUp2}, we have, for any $s \in \Gamma$:
\beq \label{BoundMonotoneCore}
\ba{rcl}
\hat{\sigma}_{p - 1}'(\|y - x\|) 
& \leq & 
\hat{\sigma}_{p - 1}'(s) 
+ \frac{ \hat{\sigma}_{p - 1}'(s)  }{ s^{p} } \|y - x\|^p \\
\\
& \overset{(\ref{Compat})}{\leq} &
\frac{(p + 1) \hat{\sigma}_{p - 1}(s)}{s}
+ \frac{(p + 1) \hat{\sigma}_{p - 1}(s)}{s^{p + 1}} \|y - x\|^p.
\ea
\eeq
Let us choose $s = \bigl[ \frac{2 \delta}{5M} \bigr]^{\frac{1}{p + 1}} \in \Gamma$.
Then, by our assumption,
\beq \label{BoundAlpha1}
\ba{rcl}
\frac{(p + 1) \hat{\sigma}_{p - 1}(s)}{s}
& = &
(p + 1)\bigl[ \frac{5M}{2 \delta} \bigr]^{\frac{1}{p + 1}}
\hat{\sigma}_{p - 1}\Bigl( \Bigl[ \frac{2\delta}{5M} \Bigr]^{\frac{1}{p + 1}} \Bigr) \\
\\
& \overset{(\ref{MBigMonotone})}{\leq} &
(p + 1)\bigl[ \frac{5M}{2 \delta} \bigr]^{\frac{1}{p + 1}} 
\frac{2\delta}{5(p + 1)} 
\;\; \overset{(\ref{AlphaChoiceMonotone})}{=} \;\; \alpha.
\ea
\eeq
Moreover,
$$
\ba{rcl}
\frac{(p + 1) \hat{\sigma}_{p - 1}(s)}{s^{p + 1}}
& \overset{(\ref{BoundAlpha1})}{\leq} &
\frac{\alpha}{s^{p}}
\;\; \equiv \;\; M.
\ea
$$
Therefore, combining the last two inequalities with~\eqref{BoundMonotoneCore},
we ensure that~\eqref{SufficientMonotone} is satisfied.
\qed

\begin{remark}
	Note that condition~\eqref{MBigMonotone}
	on our choice of $M$ is the same, up to a numerical constant in the right-hand side, 
	as condition~\eqref{MLarge} from our main Lemma~\ref{LemmaReducedGrad} on the reduced operator.
	Therefore, by ensuring condition~\eqref{MBigMonotone}, we immediately obtain both
	the monotonicity of the subproblem, and sufficient progress for one step of the method. 
\end{remark}

We conclude this section by describing a possible implementation of each iteration of algorithm~\eqref{AlgorithmMain},
in view of the result of Lemma~\ref{LemmaMonotone}.
In the most important case, $p = 1$, the subproblem in~\eqref{StatCond} 
is always monotone and can be solved by the methods of Convex Optimization.

For general $p \geq 2$ and an arbitrary value $M > 0$, 
the subproblem in~\eqref{StatCond} can be \textit{non-monotone}, even for monotone $V$, and thus
may admit multiple solutions. The standard approach is to allow inexact solutions for the subproblem
and to ensure that the inner solver admits a stopping criterion via condition~\eqref{MainProgress}
for a candidate solution $x^+$. If the complexity of the inner solver required to achieve~\eqref{MainProgress} is known, 
we can run it for a fixed budget of iterations in Step~1 of our algorithm and then check condition~\eqref{MainProgress}.
If the condition is not satisfied, the value of $M$ must be doubled and the solver relaunched.
By our theory, we ensure that the value of $M$ satisfying~\eqref{MBigMonotone} will provide the desired progress with the monotone subproblem, and it can be reached within a logarithmic number of steps of this procedure.
At this monent, the development of such schemes for $p \geq 2$ is the most challenging research direction.

\end{document}